\documentclass[12pt,twoside]{article}

\begin{document}

\centerline{}

\centerline {\Large{\bf Intuitionistic Fuzzy Banach Algebra}}

\centerline{}

\newcommand{\mvec}[1]{\mbox{\bfseries\itshape #1}}

\centerline{\bf {Bivas Dinda, T.K. Samanta and U.K. Bera }}

\centerline{}

\centerline{Department of Mathematics,}
\centerline{Mahishamuri Ramkrishna
Vidyapith, West Bengal, India. }
\centerline{e-mail: bvsdinda@gmail.com}

\centerline{Department of Mathematics, Uluberia
College, West Bengal, India.}
\centerline{e-mail: mumpu$_{-}$tapas5@yahoo.co.in}

\centerline{Department of Mathematics, City
College, Kolkata-700009, India.}
\centerline{e-mail: uttamkbera@gmail.com}
\centerline{}

\newtheorem{Theorem}{\quad Theorem}[section]

\newtheorem{definition}[Theorem]{\quad Definition}

\newtheorem{theorem}[Theorem]{\quad Theorem}

\newtheorem{remark}[Theorem]{\quad Remark}

\newtheorem{corollary}[Theorem]{\quad Corollary}

\newtheorem{note}[Theorem]{\quad Note}

\newtheorem{lemma}[Theorem]{\quad Lemma}

\newtheorem{example}[Theorem]{\quad Example}

\begin{abstract}
\textbf{\emph{Intuitionistic fuzzy Banach algebra is introduced and a few properties of it is studied. The properties of invertible elements and relation among invertible elements, open set, closed set are emphasized. Topological divisors of zero is defined and its relation with closed set are studied.}}
\end{abstract}

{\bf Keywords:}  \emph{Intuitionistic fuzzy Banach algebra, invertible elements, open set, close set, topological divisors of zero.}\\
\textbf{2010 Mathematics Subject Classification:} 46G05, 03F55.

\section{Introduction}
Fuzzy set theory is a useful tool to describe the situation in which data are imprecise or vague or uncertain. Intuitionistic fuzzy set theory handle the situation by attributing a degree of membership and a degree of non-membership to which a certain object belongs to a set. It has a wide range of application in the field of population dynamics \cite{Barros LC}, computer programming \cite{Giles}, medicine \cite{Barro} etc. \\
The concept of intuitionistic fuzzy set, as a generalisation of fuzzy sets \cite{zadeh}
was introduced by Atanassov in \cite{Atanassov}. The concept of fuzzy norm was introduced by Katsaras
\cite{Katsaras} in 1984. In 1992, Felbin\cite{Felbin1} introduced
the idea of fuzzy norm on a linear space. Cheng-Moderson
\cite{Shih-chuan} introduced another idea of fuzzy norm on a linear
space whose associated metric is same as the associated metric of
Kramosil-Michalek \cite{Kramosil}. Latter on Bag and Samanta
\cite{Bag} modified the definition of fuzzy norm on a linear space of
Cheng-Moderson \cite{Shih-chuan}.\\
In this paper, we introduce intuitionistic fuzzy Banach algebra. In section 3 we define intuitionistic fuzzy normed algebra, intuitionistic fuzzy Banach algebra. In section 4 we study the properties of invertible elements and its relation with open sets and close sets. Lastly we define topological divisor of zero and study its relation with close sets.

\section{Preliminaries}
We quote some definitions and statements of a few theorems
which will be needed in the sequel.

\begin{definition}
\cite{Schweizer} A binary operation \, $\ast \; : \; [\,0 \; , \;
1\,] \; \times \; [\,0 \; , \; 1\,] \;\, \rightarrow \;\, [\,0
\; , \; 1\,]$ \, is continuous \, $t$ - norm if \,$\ast$\, satisfies
the
following conditions \, $:$ \\
$(\,i\,)$ \hspace{0.5cm} $\ast$ \, is commutative and associative ,
\\ $(\,ii\,)$ \hspace{0.4cm} $\ast$ \, is continuous , \\
$(\,iii\,)$ \hspace{0.2cm} $a \;\ast\;1 \;\,=\;\, a \hspace{1.2cm}
\forall \;\; a \;\; \varepsilon \;\; [\,0 \;,\; 1\,]$ , \\
$(\,iv\,)$ \hspace{0.2cm} $a \;\ast\; b \;\, \leq \;\, c \;\ast\; d$
\, whenever \, $a \;\leq\; c$  ,  $b \;\leq\; d$  and  $a \,
, \, b \, , \, c \, , \, d \;\, \varepsilon \;\;[\,0 \;,\; 1\,]$.
\end{definition}
A few examples of continuous t-norm are $\,a\,\ast\,b\,=\,ab,\;\,a\,\ast\,b\,=\,min\{a,b\},\;\,a\,\ast\,b\,=\,max\{a+b-1,0\}$.

\begin{definition}
\cite{Schweizer}. A binary operation \, $\diamond \; : \; [\,0 \; ,
\; 1\,] \; \times \; [\,0 \; , \; 1\,] \;\, \rightarrow \;\,
[\,0 \; , \; 1\,]$ \, is continuous \, $t$-conorm if
\,$\diamond$\, satisfies the
following conditions \, $:$ \\
$(\,i\,)\;\;$ \hspace{0.1cm} $\diamond$ \, is commutative and
associative ,
\\ $(\,ii\,)\;$ \hspace{0.1cm} $\diamond$ \, is continuous , \\
$(\,iii\,)$ \hspace{0.1cm} $a \;\diamond\;0 \;\,=\;\, a
\hspace{1.2cm}
\forall \;\; a \;\; \in\;\; [\,0 \;,\; 1\,]$ , \\
$(\,iv\,)$ \hspace{0.1cm} $a \;\diamond\; b \;\, \leq \;\, c
\;\diamond\; d$ \, whenever \, $a \;\leq\; c$  ,  $b \;\leq\; d$
 and  $a \, , \, b \, , \, c \, , \, d \;\; \in\;\;[\,0
\;,\; 1\,].$
\end{definition}
A few examples of continuous t-conorm are $\,a\,\diamond\,b\,=\,a+b-ab,\;\,a\,\diamond\,b\,=\,max\{a,b\},\;\,a\,\diamond\,b\,=\,min\{a+b,1\}$.

\begin{definition}
\cite{Samanta} Let \,$\ast$\, be a continuous \,$t$-norm ,
\,$\diamond$\, be a continuous \,$t$-conorm  and \,$V$\, be a
linear space over the field \,$F \;(\, = \; \mathbf{R} \;\, or \;\,
\mathbf{C} \;)$. An \textbf{intuitionistic fuzzy norm} on \,$V$\,
is an object of the form  $\\A\;=\{\;(\,(\,x \;,\; t\,)
\;,\;\mu\,(\,x \;,\; t\,) \;,
\; \nu\,(\,x \;,\; t\,)\;)\;:
\;(\,x \;,\; t\,) \;\,\in\; V \;\times\;
\mathbf{R^{\,+}} \;\}$ , where $\mu \,,\, \nu\;are\; fuzzy\; sets
\;on \,$V$  \;\times\; \mathbf{R^{\,+}}$ , \,$\mu$\, denotes the
degree of membership and \,$\nu$\, denotes the degree of non -
membership \,$(\,x \;,\; t\,) \;\,\in\;\, V \;\times\;
\mathbf{R^{\,+}}$\, satisfying the following conditions $:$ \\\\
$(\,i\,)$ \hspace{0.10cm}  $\mu\,(\,x \;,\; t\,) \;+\; \nu\,(\,x
\;,\; t\,) \;\,\leq\;\, 1 \hspace{1.2cm} \forall \;\; (\,x \;,\;
t\,)
\;\,\in\;\, V \;\times\; \mathbf{R^{\,+}}\, ;$ \\
$(\,ii\,)$ \hspace{0.10cm}$\mu\,(\,x \;,\; t\,) \;\,>\;\, 0 \, ;$ \\
$(\,iii\,)$ $\mu\,(\,x \;,\; t\,) \;\,=\;\, 1\;$ if
and only if \, $x \;=\; \theta \,$, $\theta$ is null vector ; \\
$(\,iv\,)$\hspace{0.05cm} $\mu\,(\,c\,x \;,\; t\,) \;\,=\;\,
\mu\,(\,x \;,\; \frac{t}{|\,c\,|}\,)\;\;\;\;\forall\; c
\;\,\in\;\, F \, $ and $c \;\neq\; 0 \;;$ \\ $(\,v\,)$ \hspace{0.10cm} $\mu\,(\,x
\;,\; s\,) \;\ast\; \mu\,(\,y \;,\; t\,) \;\,\leq\;\, \mu\,(\,x
\;+\; y \;,\; s \;+\; t\,) \, ;$ \\ $(\,vi\,)$ \hspace{0.05cm}
$\mu\,(\,x \;,\; \cdot\,)$ is non-decreasing function of
$\;\mathbf{R^{\,+}}\,$ and  $\,\mathop {\lim }\limits_{t\;\, \to
\,\;\infty } \;\,\,\mu\,\left( {\;x\;,\;t\,} \right)=1 ;$
\\ $(\,vii\,)$ \hspace{0.10cm}$\nu\,(\,x \;,\; t\,) \;\,<\;\, 1 \, ;$ \\
$(\,viii\,)$ $\nu\,(\,x \;,\; t\,) \;\,=\;\, 0\;$  if
and only if  $\,x \;=\; \theta \, ;$ \\ $(\,ix\,)$
\hspace{0.05cm} $\nu\,(\,c\,x \;,\; t\,) \;\,=\;\, \nu\,(\,x \;,\;
\frac{t}{|\,c\,|}\,)\;\;\;\;\forall\; c
\;\,\in\;\, F \, $ and $c \;\neq\; 0 \;;$ \\ $(\,x\,)$ \hspace{0.15cm} $\nu\,(\,x
\;,\; s\,) \;\diamond\; \nu\,(\,y \;,\; t\,) \;\,\geq\;\, \nu\,(\,x
\;+\; y \;,\; s \;+\; t\,) \, ;$ \\ $(\,xi\,)$ \hspace{0.04cm}
$\nu\,(\,x \;,\; \cdot\,)$ is non-increasing function of \,
$\mathbf{R^{\,+}}\;$  and  $\;\mathop {\lim }\limits_{t\;\, \to
\,\;\infty } \;\,\,\nu\,\left( {\;x\;,\;t\,} \right)=0.$
\end{definition}

\begin{definition}
\cite{Samanta} If $A$ is an intuitionistic fuzzy norm on a linear
space $V$ then $(V\;,\;A)$ is called an intuitionistic fuzzy normed
linear space.
\end{definition}

\begin{definition}
\cite{Samanta} A sequence $\{x_n\}_n$ in an intuitionistic fuzzy normed linear space $(V\,,\,A)$ is said to \textbf{converge} to $x\;\in\;V$ if for given $r>0,\;t>0,\;0<r<1$, there exist an integer $n_0\;\in\;\mathbf{N}$ such that \\
$\;\mu\,(\,x_n\,-\,x\,,\,t\,)\;>\;1\,-\,r\;\;$
 and $\;\;\nu\,(\,x_n\,-\,x\,,\,t\,)\;<\;r\;\;$  for all $n\;\geq \;n_0$.
\end{definition}

\begin{definition}
\cite{Samanta} Let $(\;U\;,\;A\;)$ and $(\;V\;,\;B\;)$ be two
intuitionistic fuzzy normed linear space over the same field $F$. A
mapping $f$ from $(\;U\;,\;A\;)$ to $(\;V\;,\;B\;)$ is said to be \textbf{
intuitionistic fuzzy continuous} at $x_0\;\in\;U$, if for any given
$\epsilon\;>\;0\;,\alpha\;\in\;(0,1)\;,\exists\;\delta
\;=\delta(\alpha,\epsilon)\;>0\;,\beta\;=\beta(\alpha,\epsilon)\;\in\;(0,1)$
such that for all $x\;\in\;U$,
\[\mu_U(x-x_0 \;,\;\delta)\;>\;1-\beta\;\Rightarrow\;
\mu_V(f(x)-f(x_0) \;,\;\epsilon)\;>\;1-\alpha\] \[
\nu_U(x-x_0 \;,\;\delta)\;<\;\beta\;\Rightarrow\;
\nu_V(f(x)-f(x_0) \;,\;\epsilon)\;<\;\alpha\;. \]
\end{definition}

\begin{definition}
\cite{dinda} Let $\,0\;<\;r\,<\,1\;,\;t\;\in\;\mathbf{R}^+$ and $\;x\;\in\;V$.
Then the set \\
         \[\;B(\,x\,,\,r\,,\,t\,)\;=\;\{\;y\;\in\;V\;:\;\mu(\,x\,-\,y \,,\, t\,)
         \;>\;1-r\;,\;\;\nu(\,x-y \,,\, t\,)\;<\;r\;\}\]
is called an \textbf{open ball} in \,$(\,V\;,\;A\,)$ \,with $x$ as its
center and $r$ as its radious with respect to $\,t\,$.
\end{definition}

\section{Intuitionistic Fuzzy Banach Algebra}
\begin{definition}
 Let \,$\ast$\, be a continuous \,$t$-norm ,
\,$\diamond$\, be a continuous \,$t$-conorm  and \,$\mathcal{A}$\, be an
algebra over the field \,$k\;(\, = \; \mathbf{R} \;\, or \;\,
\mathbf{C} \;)$. An \textbf{intuitionistic fuzzy normed algebra} is an object of the form  $\left(\,\mathcal{A},\,\mu,\,\nu,\,\ast,\,\diamond\right)$ , where $\mu \,,\, \nu\;$ are fuzzy sets
on \,$V \;\times\; \mathbf{R^{\,+}}$, \,$\mu$\, denotes the
degree of membership and $\nu$\, denotes the degree of non-membership satisfying the following conditions for every $\;x,y\,\in\,\mathcal{A}$ and $s,t\,\in\,\mathbf{R^+}$ ;

$(i)$ \hspace{0.10cm}  $\mu\,(\,x \;,\; t\,) \;+\; \nu\,(\,x
\;,\; t\,) \;\,\leq\;\, 1 $

$(ii)$ \hspace{0.10cm}$\mu\,(\,x \;,\; t\,) \;\,>\;\, 0 \, ;$

$(iii)$ $\mu\,(\,x \;,\; t\,) \;\,=\;\, 1\;$ if
and only if \, $x \;=\; \theta \,$, $\theta$ is null vector;

$(iv)$\hspace{0.05cm} $\mu\,(\,c\,x \;,\; t\,) \;\,=\;\,
\mu\,(\,x , \frac{t}{|\,c\,|}\,)\;\;\;\;\forall\; c
\;\,\in\;\, F \, $ and $c \;\neq\; 0 ;$

$(v)$ \hspace{0.10cm} $\mu\,(\,x
\;,\; s\,) \,\ast\, \mu\,(\,y \;,\; t\,) \;\,\leq\;\, \mu\,(\,x
\;+\; y \;,\; s \;+\; t\,) \, ;$

$(vi)$ $\,\max\{\mu(x,s),\mu(y,t)\}\,\leq\,\mu(xy,s+t);$

$(vii)$ $\,\mathop {\lim }\limits_{t\;\, \to
\,\;\infty } \;\,\,\mu\,\left( {\;x\;,\;t\,} \right)=1\,$ and $\,\mathop {\lim }\limits_{t\;\, \to
\,\;0} \;\,\,\mu\,\left( {\;x\;,\;t\,} \right)=0 $

$(\,viii\,)$ \hspace{0.10cm}$\nu\,(\,x \;,\; t\,) \;\,<\;\, 1 \, ;$

$(ix)$ $\nu\,(\,x \;,\; t\,) \;\,=\;\, 0\;$  if
and only if  $\,x \;=\; \theta ,\;\theta$ is null vector;

$(x)$\hspace{0.05cm} $\nu\,(\,c\,x \;,\; t\,) \;\,=\;\, \nu\,(\,x \;,\;
\frac{t}{|\,c\,|}\,)\;\;\;\;\forall\; c
\;\,\in\;\, F \, $ and $c \;\neq\; 0 \;;$

$(xi)$ \hspace{0.15cm} $\nu\,(\,x
\;,\; s\,) \;\diamond\; \nu\,(\,y \;,\; t\,) \;\,\geq\;\, \nu\,(\,x
\;+\; y \;,\; s \;+\; t\,) ;$

$(xii)$ $\,\min\{\nu(x,s),\nu(y,t)\}\,\geq\,\nu(xy,s+t);$

$(xiii)$  $\;\mathop {\lim }\limits_{t\;\, \to
\,\;\infty } \;\,\,\nu\,\left( {\;x\;,\;t\,} \right)=0\;$ and $\;\mathop {\lim }\limits_{t\;\, \to
\,\;0 } \;\,\,\nu\,\left( {\;x\;,\;t\,} \right)=1\,.$
\end{definition}

For an intuitionistic fuzzy normed algebra $\;\left(\,\mathcal{A},\,\mu,\,\nu,\,\ast,\,\diamond\right)\,$ we further assume that:
$\;\;\;(xiv)\;\left. {{}_{a\;\; \ast \;\;a\;\; =
\;\;a}^{a\;\; \diamond \;\;a\;\; = \;\;a} \;\;}
\right\}\;\;\;\,\;\;\;\;\,,\;\;\,$ for all
$\;a\;\; \in \;\;[\,0\;\,,\;\,1\,].$

\begin{definition}
A sequence $\{x_n\}_n$ in an intuitionistic fuzzy normed algebra $\left(\,\mathcal{A},\,\mu,\,\nu,\,\ast,\,\diamond\right)$ is said to \textbf{converge} to $\,x\;\in\,\mathcal{A}\,$ if for given $r>0,\;t>0,\;0<r<1$, there exist an integer $n_0\;\in\;\mathbf{N}$ such that \\
$\;\mu\,(\,x_n\,-\,x\,,\,t\,)\;>\;1\,-\,r\;\;$
 and $\;\;\nu\,(\,x_n\,-\,x\,,\,t\,)\;<\;r\;\;$  for all $n\;\geq \;n_0$.
\end{definition}

\begin{theorem}
In an intuitionistic fuzzy normed algebra\,$\left(\,\mathcal{A},\,\mu,\,\nu,\,\ast,\,\diamond\right)$, a sequence \,$\left\{ {\;x_{\,n}
\;} \right\}_{\,n} $\, converges to \,$x \;\,\varepsilon\;\,\mathcal{ A }$\,
if and only if \,$\mathop {\lim }\limits_{n\;\, \to \,\;\infty }
\,\mu\,\left( {\;x_{\,n} \; - \;x\;,\;t\,} \right)\;= \;1$\, and
\,$\mathop {\lim }\limits_{n\;\, \to \,\;\infty } \,\nu\,\left(
{\;x_{\,n} \; - \;x\;,\;t\,} \right)\; =\;0$ .
\end{theorem}

\begin{definition}
 A sequence $\{x_n\}_n$ in an intuitionistic fuzzy normed algebra $\left(\,\mathcal{A},\,\mu,\,\nu,\,\ast,\,\diamond\right)$ is said to be \textbf{cauchy sequence} if $\mathop {\lim }\limits_{n\;\,
\to \,\;\infty } \,\,\mu(x_{n+p}-x_n ,t)=1\; $ and $\mathop {\lim }\limits_{n\;\,
\to \,\;\infty } \,\,\nu(x_{n+p}-x_n ,t)=0\;\;,\,t\,\in\,\mathbf{R^+},\;p=1,2,3,..... $.
\end{definition}

\begin{definition}
An intuitionistic fuzzy normed algebra $\left(\,\mathcal{A},\,\mu,\,\nu,\,\ast,\,\diamond\right)$ is said to be \textbf{complete} if every cauchy sequence in $\mathcal{A}$ converges to an element of $\mathcal{A}$.
\end{definition}

\begin{definition}
A complete intuitionistic fuzzy normed algebra is called intuitionistic fuzzy Banach algebra.
\end{definition}

\begin{theorem}
In an intuitionostic fuzzy Banach algebra $\,(\mathcal{A},\,\mu,\,\nu,\,\ast,\,\diamond)\,$ two sequences $\,\{x_n\}_n\,$ and $\,\{y_n\}_n\,$ be such that $\,x_n\,\rightarrow\,x\,$ and $\,y_n\,\rightarrow\,y\,$ then $\,x_ny_n\,\rightarrow\,xy\,.$
\end{theorem}
{\bf Proof.} For any $\,t>0\,,\,$ $x_n\rightarrow\,x\,$ implies
\[\mathop {\lim }\limits_{n\;\to \,\;\infty }\,\mu\,\left( {\;x_{\,n} \; - \;x\;,\;t\,} \right)\,= \,1\,,\;\;\;
\,\mathop {\lim }\limits_{n\;\, \to \,\;\infty } \,\nu\,\left(
{\;x_{\,n} \; - \;x\;,\;t\,} \right)\; =\;0\]
and $\,y_n\rightarrow\,y\,$ implies
\[\mathop {\lim }\limits_{n\;\to \,\;\infty }\,\mu\,\left( {\;y_{\,n} \; - \;y\;,\;t\,} \right)\,= \,1\,,\;\;\;
\,\mathop {\lim }\limits_{n\;\, \to \,\;\infty } \,\nu\,\left(
{\;y_{\,n} \; - \;y\;,\;t\,} \right)\; =\;0\]
Now \[\mathop {\lim }\limits_{n\;\to \,\;\infty }\,\mu\,\left( \,x_n\,y_n \, - \,xy\,,\,t\, \right)\,\geq\,\mathop {\lim }\limits_{n\to \,\infty }\,\mu\,\left( \,(x_{\,n} \; - \;x)y_n\,,\,\frac{t}{2}\, \right)\,\ast\,\mathop {\lim }\limits_{n\;\to \,\;\infty }\,\mu\,\left(\, (y_{\,n} \; - \;y)x\,,\,\frac{t}{2}\, \right)\]
\[\geq\,\max\;\left\{\mathop {\lim }\limits_{n\to \,\infty }\,\mu\,\left( \,x_{\,n} \; - \;x\,,\,\frac{t}{4}\, \right),\,\mathop {\lim }\limits_{n\to \,\infty }\,\mu\,\left( \,y_n\,,\,\frac{t}{4}\, \right)\right\}\hspace{7 cm}\]
\[\,\ast\,\max\;\left\{\mathop {\lim }\limits_{n\to \,\infty }\,\mu\,\left( \,y_{\,n} \; - \;y\,,\,\frac{t}{4}\, \right),\,\mathop {\lim }\limits_{n\to \,\infty }\,\mu\,\left( \,x\,,\,\frac{t}{4}\, \right)\right\}\hspace{-2 cm}\]
\[=\,\max\;\left\{\,1\,,\,\mathop {\lim }\limits_{n\to \,\infty }\,\mu\,\left( \,y_n\,,\,\frac{t}{4}\, \right)\right\}\,\ast\,\max\;\left\{\,1\,,\,\mathop {\lim }\limits_{n\to \,\infty }\,\mu\,\left( \,x\,,\,\frac{t}{4}\, \right)\right\}\hspace{7 cm}\]
$=\,1\,\ast\,1\;=\,1\,.$\\
and
\[\mathop {\lim }\limits_{n\;\to \,\;\infty }\,\nu\,\left( \,x_n\,y_n \, - \,xy\,,\,t\, \right)\,\leq\,\mathop {\lim }\limits_{n\to \,\infty }\,\nu\,\left( \,(x_{\,n} \; - \;x)y_n\,,\,\frac{t}{2}\, \right)\,\diamond\,\mathop {\lim }\limits_{n\;\to \,\;\infty }\,\nu\,\left(\, (y_{\,n} \; - \;y)x\,,\,\frac{t}{2}\, \right)\]
\[\leq\,\min\;\left\{\mathop {\lim }\limits_{n\to \,\infty }\,\nu\,\left( \,x_{\,n} \; - \;x\,,\,\frac{t}{4}\, \right),\,\mathop {\lim }\limits_{n\to \,\infty }\,\nu\,\left( \,y_n\,,\,\frac{t}{4}\, \right)\right\}\hspace{7 cm}\]
\[\,\diamond\,\min\;\left\{\mathop {\lim }\limits_{n\to \,\infty }\,\nu\,\left( \,y_{\,n} \; - \;y\,,\,\frac{t}{4}\, \right),\,\mathop {\lim }\limits_{n\to \,\infty }\,\nu\,\left( \,x\,,\,\frac{t}{4}\, \right)\right\}\hspace{-2 cm}\]
\[=\,\min\;\left\{\,0\,,\,\mathop {\lim }\limits_{n\to \,\infty }\,\nu\,\left( \,y_n\,,\,\frac{t}{4}\, \right)\right\}\,\diamond\,\min\;\left\{\,0\,,\,\mathop {\lim }\limits_{n\to \,\infty }\,\nu\,\left( \,x\,,\,\frac{t}{4}\, \right)\right\}\hspace{7 cm}\]
$=\,0\,\diamond\,0\;=\,0\,.$\\
This complites the proof.

\section{Invertible elements and Topological divisors of Zero}
\begin{theorem}
Let $\,\left(\mathcal{A},\,\mu,\,\nu,\,\ast,\,\diamond\right)\,$ be an intuitionistic fuzzy Banach algebra. If $\,x\,\in\,\mathcal{A}\,$ be such that $\,x\,\in\,B(\theta,r,t),\,0<r<1\,$ then $\,(e-x)\,$ is invertible and $\,(e-x)^{-1}\,=\,e+\sum_{{_(n=1)}}^{^\infty} x^n $
\end{theorem}
{\bf Proof.} Since $\,x\,\in\,B(\theta,r,t)\,$ we have for any $\,t\,\in\,\mathbf{R}^+, \;\\ \mu(x,t)>1-r,\;\nu(x,t)<r\\$
To prove the theorem we first show that the series $\,\sum_{{_(n=1)}}^{^\infty} x^n\,$ is convergent to some element of $\,\mathcal{A}\,.\,$ Let $\,s_n\,=\sum_{{_(k=1)}}^{^n} x^k,\,$ then it is sufficient to prove that \[\mu\left(s_{n+p}-s_n,t\right)>1-r\;\;\;and\;\;\;\nu\left(s_{n+p}-s_n,t\right)<r\]
Now \[\mu\left(s_{n+p}-s_n,t\right)\,=\,\mu\left(x^{n+1}+x^{n+2}+\cdots+x^{n+p},t\right)\hspace{5 cm}\]
\[\geq\,\mu(x^{n+1},t_1)\,\ast\,\mu(x^{n+2},t_2)\,\ast\,\cdots\,\mu(x^{n+p},t_p)\,,\;\;\;\;where \;\;\;t_1+t_2+\cdots+t_p\,=\,t\]
\[\geq\;\max\,\left\{\,\mu(x,t_{11}),\,\mu(x,t_{12}),\,\cdots\,\mu(x,t_{1n+1})\,\right\}\,\ast\,
\max\,\left\{\,\mu(x,t_{21}),\,\mu(x,t_{22}),\,\cdots\,\mu(x,t_{2n+2})\,\right\}\,\]\[\ast\,\cdots\cdots\,\ast\,
\max\,\left\{\,\mu(x,t_{p1}),\,\mu(x,t_{p2}),\,\cdots\,\mu(x,t_{pn+p})\,\right\}\,,\]
\[\;where \;\sum_{{_(j=1)}}^{^n+i} t_{ij}\,=\,t_i\,,\,i=1,2,\cdots,p\,\hspace{-6 cm}\]
\[>\,(1-r)\,\ast\,(1-r)\,\ast\,\cdots\,\ast\,(1-r)\,=\,(1-r)\hspace{8 cm}\]
and
\[\nu\left(s_{n+p}-s_n,t\right)\,=\,\nu\left(x^{n+1}+x^{n+2}+\cdots+x^{n+p},t\right)\hspace{5 cm}\]
\[\leq\,\nu(x^{n+1},t_1)\,\diamond\,\nu(x^{n+2},t_2)\,\diamond\,\cdots\,\nu(x^{n+p},t_p)\,,\;\;\;\;where \;\;\;t_1+t_2+\cdots+t_p\,=\,t\]
\[\leq\;\min\,\left\{\,\nu(x,t_{11}),\,\nu(x,t_{12}),\,\cdots\,\nu(x,t_{1n+1})\,\right\}\,\diamond\,
\min\,\left\{\,\nu(x,t_{21}),\,\nu(x,t_{22}),\,\cdots\,\nu(x,t_{2n+2})\,\right\}\,\]\[\diamond\,\cdots\cdots\,\diamond\,
\min\,\left\{\,\nu(x,t_{p1}),\,\nu(x,t_{p2}),\,\cdots\,\nu(x,t_{pn+p})\,\right\}\,,\]
\[\;where \;\sum_{{_(j=1)}}^{^n+i} t_{ij}\,=\,t_i\,,\,i=1,2,\cdots,p\,\hspace{-6 cm}\]
\[<\,r\,\diamond\,r\,\diamond\,\cdots\diamond r\,=\,r\hspace{10 cm}\]
Thus the series $\,\sum_{{_(n=1)}}^{^\infty} x^n\,$ is convergent. Since $\,\mathcal{A}\,$ is complete, $\,\sum_{{_(n=1)}}^{^\infty} x^n\,$ converges to some element of $\,\mathcal{A}\,.\\$
Let $\,s\,=\,e\,+\,\sum_{{_(n=1)}}^{^\infty} x^n\,.\\$
Now $\,(e-x)(e+x+\cdots+x^n)\,=\,(e+x+\cdots+x^n)(e-x)\,=\,e-x^{n+1}\,.\\$
Also, $\,\mu\left(x^{n+1},t\right)\,\geq \;\max\left\{\mu(x,t_1),\,\mu(x,t_2),\,\cdots,\,\mu(x,t_{n+1})\,\right\}\,>\,1-r\,,\;$ and $\;\nu\left(x^{n+1},t\right)\,\leq \;\min\left\{\nu(x,t_1),\,\nu(x,t_2),\,\cdots,\,\nu(x,t_{n+1})\,\right\}\,<\,r\,,\;\;\;where \;\;t_1+t_2+\cdots+t_{n+1}\,=\,t\,.\\$
So, $\,x^{n+1}\,\rightarrow\,\theta\;\;as \;n\rightarrow\,\infty\,.\;$ Therefore letting $\;n\rightarrow\,\infty\,$ and remembering that multiplication on $\,\mathcal{A}\,$ is continuous, we obtain $\;(e-x)s\,=\,s(e-x)\,=\,e\,.\\$
Hence $\,s\,=\,(e-x)^{-1}\,.\;$ This complites the proof.

\begin{corollary}\label{c1}
Let $\,x\,\in\,\mathcal{A}\,$ be such that $\,e-x\in B(\theta,r,t),\,0<r<1\,$ then $\,x^{-1}\,$ exist and $\,x^{-1}\,=\,e+\sum_{{_(n=1)}}^{^\infty}(e-x)^n\,.$
\end{corollary}

\begin{corollary}
Let $\,x\in\,\mathcal{A}\,$ and $\,\lambda(\neq\,0)\,$ be a scaler such that $\,\frac{x}{\lambda}\,\in\,B(\theta,r,t),\,0<r<1.\,$ Then $\,(\lambda\,e-x)^{-1}\,$ exists and $\,(\lambda\,e-x)^{-1}\,=\,\sum_{{_(n=1)}}^{^\infty} \lambda^{-n}x^{n-1}$
\end{corollary}
{\bf Proof.} Since $\,\frac{x}{\lambda}\in\,B(\theta,r,t)\,,\,$ for any $\,t\in \,\mathbf{R}^+\,$ we have $\,\mu(\frac{x}{\lambda},t)>1-r\;\;and\;\;\nu(\frac{x}{\lambda},t)<r\,.\\$
If $\,y\in\,\mathcal{A}\,$ be such that $\,y^{-1}\in\,\mathcal{A}\,$ and $\,c(\neq\,0)\,$ be a scaler then it is clear that $\,(cy)^{-1}\,$ exist and $\,(cy)^{-1}\,=\,c^{-1}y^{-1}\,.\\$
Now $\,(\lambda\,e-x)\,=\,\lambda(e-\frac{x}{\lambda})\,$ and we show that $\,(e-\frac{x}{\lambda})^{-1}\,$ exists.
$\,\mu\left(e-(e-\frac{x}{\lambda}),t\right)\,=\,\mu(\frac{x}{\lambda},t)\,>1-r\;\;and\;\;\nu\left(e-(e-\frac{x}{\lambda}),t\right)\,=\, \nu(\frac{x}{\lambda},t)\, <r\,,\;$ by hypothesis.\\
So, by corollary \ref{c1} $\,(e-\frac{x}{\lambda})^{-1}\,$ exists and so $\,(\lambda\,e-x)^{-1}\,$ exists.
\[\left(\,\lambda\,e-x\,\right)^{-1}\,=\,\left\{\,\lambda(e-\frac{x}{\lambda})\,\right\}\,=\,\lambda^{-1}\,\left(e-\lambda^{-1}x\right)^{-1}\hspace{8 cm}\]
\[=\,\lambda^{-1}\left(e\,+\,\sum_{{_(n=1)}}^{^\infty}[e-(e-\lambda^{-1}x)]^n\right)\hspace{2.75 cm}\]
\[=\,\lambda^{-1}\left(e\,+\,\sum_{{_(n=1)}}^{^\infty}\,\lambda^{-n}x^n\right)\,=\,\sum_{{_(n=1)}}^{^\infty}\,\lambda^{-n}x^{n-1}\hspace{1.5 cm}\]
This complites the proof.

\begin{theorem}
The set of all invertible elements of an intuitionistic fuzzy Banach algebra $\,(\mathcal{A},\mu,\nu,\ast,\diamond)\,$ is an open subset of $\,\mathcal{A}.$
\end{theorem}
{\bf Proof.} Let G be the set of all invertible elements of $\mathcal{A}.\;$ Let $\,x_0\in\,G\,,\,$ we have to show that the open ball with centre at $\,x_0\,$ and radious $\,r\,(0<r<1)\,$ contain in G  i.e., every point $x$ of the open ball $B(x_0,r,t)\,$ satisfies the inequality
\begin{equation}\label{e1}
\mu(x-x_0,t)>1-r\;\;,\;\;\nu(x-x_0,t)<r.
\end{equation}
Let us choose $r$ such that \[r<\min\{\mu(x_0^{-1},t),\,1-\nu(x_0^{-1},t)\,\},\;\;\;\;0<\mu(x_0^{-1},t),\,\nu(x_0^{-1},t)<1\]
Then $\,\mu(x-x_0,t)\,>\,1-\left(1-\mu(x_0^{-1},t)\right)\,=\,\mu(x_0^{-1},t)\;,\;\;\nu(x-x_0,t)\,<\,\nu(x_0^{-1},t)\,.\\$
Let $\,y\,=\,x_0 ^{-1}x\;and \;z\,=\,e-y\,.$  Then
\[\mu(z,t)\,=\,\mu(e-y,t)\,=\,\mu(y-e,t)\,=\,\mu(x_0 ^{-1}x\,-\,x_0 ^{-1}x_0,t)\,=\,\mu\left(x_0 ^{-1}(x-x_0),t\right)\]
\[\geq\,\max\,\{\mu(x_0 ^{-1},\frac{t}{2}),\,\mu(x-x_0,\frac{t}{2})\}\,=\,\mu(x-x_0,\frac{t}{2})\,>\,1-r. \hspace{6 cm}\]
\[\nu(z,t)\,=\,\nu(e-y,t)\,=\,\nu(y-e,t)\,=\,\nu(x_0 ^{-1}x\,-\,x_0 ^{-1}x_0,t)\,=\,\nu\left(x_0 ^{-1}(x-x_0),t\right)\]
\[\leq\,\min\,\{\nu(x_0 ^{-1},\frac{t}{2}),\,\nu(x-x_0,\frac{t}{2})\}\,=\,\nu(x-x_0,\frac{t}{2})\,<\,r.\hspace{6 cm}\]
Thus $\,z\in\,B(\theta,r,t)\,$ and hence $\,e-z\,$ is invertible. That is, $y$ is invertible. So, $y\in\,G.\\$
Now $\,x-\theta\in\,G\,$ and $y\in\,G$ and by our earlier verification $G$ is a group. So, $x_0 y\,=\,x_0 x_0 ^{-1}x\,=\,x\,.$ So, any point $x$ satisfies (\ref{e1}) belongs to G; that is, $B(x_0,r,t)\subseteq\,G\,.$ Hence, $G$ is an open subset of $\,\mathcal{A}.$

\begin{corollary}
The set of all non-invertible elements of an intuitionistic fuzzy Banach algebra $\,(\mathcal{A},\mu,\nu,\ast,\diamond)\,$ is a closed subset of $\,\mathcal{A}.$
\end{corollary}

\begin{theorem}
Let $\,G\,$ be the set of all invertible elements of $\,(\mathcal{A},\mu,\nu,\ast,\diamond)\,.\;$ The mapping $\,x\rightarrow\,x^{-1}\,$ of $G$ into $G$ is strongly intuitionistic fuzzy continuous.
\end{theorem}
{\bf Proof.} Let $x_0$ be an element of G. Let $x\in G$ be such that $x\in\,B(x_0,r,t),$ where $\,r\,<\,\min\{\,1-\mu(x_0^{-1},t),\;\nu(x_0^{-1},t)\,\}\,.\,$ Therefore,\\
$\mu(x-x_0,t)\,>\,1-r\,>\mu(x_0 ^{-1},t)\;\;and \;\;\nu(x-x_0,t)\,<\,r\,<\,\nu(x_0 ^{-1},t)\\$
Since $x_0\in G$ is arbitrary $\,\forall x\in G,\; $ we have\\
\[\mu(x-x_0,t)\,>\,\mu(x^{-1},t)\;\;and \;\;\nu(x-x_0,t)\,>\,\nu(x^{-1},t)\]
Let $x\in\,G\;$ then $T(x)\,=\,x^{-1}\,.$ Now\\
\[\mu(T(x)-T(x_0),\epsilon)=\mu(x^{-1}-{x_0}^{-1},\epsilon)=\mu\left(\,(x^{-1}x_0-e){x_0}^{-1},\epsilon\,\right)\]
\[\geq \;\max\{\mu(x^{-1}x_0-e,\frac{\epsilon}{2})\,,\,\mu({x_0}^{-1},\frac{\epsilon}{2})\}\]
\[\geq \;\max\{\mu\left(x^{-1}(x_0-x),\frac{\epsilon}{2}\right)\,,\,\mu({x_0}^{-1},\frac{\epsilon}{2})\}\hspace{-0.5 cm}\]
\[\geq \;\max\left\{\max\{\mu(x^{-1},\frac{\epsilon}{4}),\;\mu(x_0-x,\frac{\epsilon}{4})\}\,,\,\mu({x_0}^{-1},\frac{\epsilon}{2})\right\}\hspace{-2.5 cm}\]
\[=\max\{\mu(x-x_0,\frac{\epsilon}{4}),\mu({x_0}^{-1},\frac{\epsilon}{2})\}\hspace{1.5 cm}\]
\[=\mu(x-x_0,\frac{\epsilon}{4})\hspace{4.5 cm}\]

\[\nu(T(x)-T(x_0),\epsilon)=\nu(x^{-1}-{x_0}^{-1},\epsilon)=\nu\left(\,(x^{-1}x_0-e){x_0}^{-1},\epsilon\,\right)\]
\[\leq \;\min\{\nu(x^{-1}x_0-e,\frac{\epsilon}{2})\,,\,\nu({x_0}^{-1},\frac{\epsilon}{2})\}\]
\[\leq \;\min\{\nu\left(x^{-1}(x_0-x),\frac{\epsilon}{2}\right)\,,\,\nu({x_0}^{-1},\frac{\epsilon}{2})\}\hspace{-0.5 cm}\]
\[\leq \;\min\left\{\min\{\nu(x^{-1},\frac{\epsilon}{4}),\;\nu(x_0-x,\frac{\epsilon}{4})\}\,,\,\nu({x_0}^{-1},\frac{\epsilon}{2})\right\}\hspace{-2.5 cm}\]
\[=\min\{\nu(x-x_0,\frac{\epsilon}{4}),\nu({x_0}^{-1},\frac{\epsilon}{2})\}\hspace{1.5 cm}\]
\[=\nu(x-x_0,\frac{\epsilon}{4})\hspace{4.5 cm}\]
Hence $T$ is strongly intuitionistic fuzzy continuous at $x_0\in G.\,$ Since $x_0$ is arbitrary, $T$ is strongly intuitionistic fuzzy continuous on $G$.

\begin{definition}
Let $\,(\mathcal{A},\mu,\nu,\ast,\diamond)\,$ be an intuitionistic fuzzy Banach algebra. An element $\,z\in\,\mathcal{A}\,$ is called a {\bf topological divisor of zero} if there exist a sequence $\,\{z_n\}_n,\,z_n\in\,\mathcal{A}\,$ satisfies $\,z_n\,\not \in\,B(\theta,r,t),\,0<r<1\,\;(\,i.e., \, \mu(z_n,t)<1-r \;or\; \nu(z_n,t)>r \; or \;both \,)\,$ be such that,  {\bf either} \[\,\mathop {\lim }\limits_{n\;\to \,\;\infty }\,\mu\,\left( \,z_n\,z\,,\,t\, \right)\,=\,1\,\; and \;\,\mathop {\lim }\limits_{n\;\to \,\;\infty }\,\nu\,\left( \,z_n\,z\,,\,t\, \right)\,=\,0\,\] {\bf or} \[\,\mathop {\lim }\limits_{n\;\to \,\;\infty }\,\mu\,\left( \,z\,z_n\,,\,t\, \right)\,=\,1\,\; and \;\,\mathop {\lim }\limits_{n\;\to \,\;\infty }\,\nu\,\left( \,z\,z_n\,,\,t\, \right)\,=\,0\,.\]
\end{definition}

\begin{theorem}
Let $\,Z\,$ be the set of all topological divisors of zero in $\,\mathcal{A}\,.\;$ Then $\;Z\,\subseteq\,S\,,$ where $S$ the set of non-invertible elements of $\mathcal{A}\,.$
\end{theorem}
{\bf Proof.} Let $\,z\in\,Z\,$ then there exist a sequence $\,\{z_n\}_n,\,z_n\in\,\mathcal{A}\,$ satisfies $\,z_n\,\not \in\,B(\theta,r,t),\,0<r<1\,\;(\,i.e., \, \mu(z_n,t)<1-r \;or\; \nu(z_n,t)>r \; or \;both \,)\,$ be such that,  {\bf either} \[\,\mathop {\lim }\limits_{n\;\to \,\;\infty }\,\mu\,\left( \,z_n\,z\,,\,t\, \right)\,=\,1\,\; and \;\,\mathop {\lim }\limits_{n\;\to \,\;\infty }\,\nu\,\left( \,z_n\,z\,,\,t\, \right)\,=\,0\,\] {\bf or} \[\,\mathop {\lim }\limits_{n\;\to \,\;\infty }\,\mu\,\left( \,z\,z_n\,,\,t\, \right)\,=\,1\,\; and \;\,\mathop {\lim }\limits_{n\;\to \,\;\infty }\,\nu\,\left( \,z\,z_n\,,\,t\, \right)\,=\,0\,.\]
Suppose $\mathop {\lim }\limits_{n\;\to \,\;\infty }\,\mu\,\left( \,z\,z_n\,,\,t\, \right)\,=\,1\,\; and \;\,\mathop {\lim }\limits_{n\;\to \,\;\infty }\,\nu\,\left( \,z\,z_n\,,\,t\, \right)\,=\,0\,.\\$
Let $G$ be the set of non-invertible elements of $\mathcal{A}\,.$ If possible let $\,z\in\,G\,.$ Then $\,z^{-1}\,$ exists. Now since multiplicatin is a continuous operation, we should have \[z_n\,=\,(z^{-1}z)z_n\,=\,z^{-1} (zz_n)\rightarrow\,z^{-1}.\theta\,=\,0\;\;\; as\;\; n\rightarrow\,\infty\,.\]
which contradicts the fact that $\,z_n\,\not \in\,B(\theta,r,t),\,0<r<1\,\;(\,i.e., \, \mu(z_n,t)<1-r \;or\; \nu(z_n,t)>r \; or \;both \,)\,.\\$
So, $\,z\in\,S\,$ and this complites the proof.

\end{document}